\documentclass[11pt]{amsart}
\input xy
\xyoption{all}

\begin{document}

\newtheorem{theorem}{Theorem}[section]
\newtheorem{lemma}[theorem]{Lemma}
\newtheorem{cor}[theorem]{Corollary}
\newtheorem{prop}[theorem]{Proposition}

\newcommand{\zz}{{\mathbb Z}}
\newcommand{\zq}{{\mathbb Q}}
\newcommand{\lra}{\longrightarrow}
\newcommand{\ra}{\rightarrow}
\newcommand{\slna}{SL_n(A)}
\newcommand{\slnk}{SL_n(k)}
\newcommand{\cG}{{\mathcal G}}
\newcommand{\cT}{{\mathcal T}}
\newcommand{\cU}{{\mathcal U}}
\newcommand{\ahat}{\widehat{A}}
\newcommand{\slnahat}{SL_n(\ahat)}
\newcommand{\mx}{{\mathfrak m}_x}

\title{Relative completions and $K_2$ of curves}
\author{Kevin P.~Knudson}\thanks{Partially supported by the
National Security Agency}
\address{Department of Mathematics \& Statistics, Mississippi
State University, P.O.~Box MA, Mississippi State, MS 39762}
\email{knudson@math.msstate.edu}



\begin{abstract} We compute the completion of the special linear
group over the coordinate ring of a curve over a number field $k$
relative to its representation in $\slnk$, and relate this to the
study of $K_2$ of the curve.
\end{abstract}

\maketitle

\section*{Introduction}
The algebraic $K$-groups of curves defined over number fields are
mysterious objects and have received a great deal of scrutiny.
There is a plethora of conjectures about the groups $K_1$ and
$K_2$ and a few results supporting these guesses.  For example, if
$A$ is the coordinate ring of a smooth affine curve over a number
field $k$, then a conjecture of Vaserstein asserts that $SK_1(A)$
is torsion.  This is supported by calculations of Raskind
\cite{raskind}, who showed that $SK_1(A)\otimes \zq/\zz=0$, and by
more recent work of {\O}stvaer and Rosenschon \cite{rosen}.

The group $K_2$ is even more troublesome, but the work of several
authors (\cite{beil}, \cite{bloch}, \cite{blochgrayson},
\cite{kk2}, \cite{rs}) has shed some light on the structure of
$K_2(A)$. For example, the existence of regulator maps shows that
nontrivial elements exist; the rank of $K_2(A)$ is conjectured to
be related to the number of infinite places of $k$ when $C$ is an
elliptic curve; the second level of the rank filtration,
$r_2K_2(A)_{\zq}$, vanishes when $C$ is an elliptic curve. Still,
it is not known for {\em any} curve whether or not $K_2(A)$ has
finite rank.

In this note, we use Deligne's notion of relative completion to
study the group $K_2(A)$.  A full summary of this construction is
given in Section \ref{relative}, but the basic idea is the
following.  Let $S$ be a reductive algebraic group over a field
$F$ and let $\rho:\Gamma\ra S$ be a map of a discrete group
$\Gamma$ into $S$ with Zariski dense image.  The completion of
$\Gamma$ relative to $\rho$ is a proalgebraic group $\cG$ over $F$
together with a lift $\widetilde{\rho}:\Gamma\ra\cG$ such that the
diagram
$$\xymatrix{
\Gamma \ar[r]^{\widetilde{\rho}}\ar[dr]_{\rho} & \cG\ar[d] \\
             & S}$$
commutes.  Using a suitably defined notion of the continuous
cohomology of $\cG$, $H^\bullet_{\text{cts}}(\cG;F)$, there is an
injective map
$$H^2_{\text{cts}}(\cG;F)\lra H^2(\Gamma;F).$$  This allows one to
obtain a lower bound on the rank of $H^2(\Gamma;F)$.

Here, we consider the group $\Gamma=\slna$, $n\ge 3$, and
$\rho:\slna\ra\slnk$ given by reduction modulo the maximal ideal
of a $k$-rational point on the curve (we assume one exists, for
simplicity).  The main result of this paper is the following
calculation of the completion of $\slna$ relative to $\rho$.

\medskip

\noindent {\bf Theorem \ref{mainthm}.}  {\em The completion of
$\slna$ relative to $\rho$ is isomorphic to the group
$SL_n(k[[T]])$.}

\medskip

Unfortunately, this provides no information about the rank of
$H^2(\slna)$ (which is related to $K_2(A)$) in light of the
following result.

\medskip

\noindent {\bf Proposition \ref{h2comp}.} {\em For $n\ge 3$,
$H^2_{\text{\em cts}}(SL_n(k[[T]]);k)=0$.}

\medskip

Thus, to obtain a complete understanding of $K_2(A)$, one needs to
compute the group $H_2(\slna;\zz)$, $n\ge 3$.  such a calculation
remains elusive, however.

\medskip

\noindent {\em Acknowledgements.} The author thanks Chuck Weibel
for many useful conversations.

\medskip

\noindent {\em Notation.} Throughout this paper, $k$ denotes a
number field and $C$ denotes a smooth affine curve over $k$.  We
denote the coordinate ring $k[C]$ by $A$.  We assume, for
simplicity, that $C$ has a $k$-rational point $x$ and we denote
the associated maximal ideal of $A$ by $\mx$.  We also assume that
$SK_1(A)\otimes \zq=0$.  For a group $G$, $\Gamma^iG$ denotes the
$i$-th term of the lower central series.

\medskip

\noindent {\em Remark.} The results of this paper are also valid
for finite fields, but the interest in that case is subsumed by
Quillen's calculations \cite{quillen}.

\section{Relative Completions}\label{relative}
Relative completion is a generalization of the classical Malcev
(or unipotent) completion.  Proofs of the results in this section
may be found in \cite{hain} or \cite{knudson}.

Let $\Gamma$ be a discrete group and let $S$ be a reductive linear
algebraic group over a field $F$.  Suppose that $\rho:\Gamma\ra S$
is a homomorphism with Zariski dense image.  The {\em completion
of $\Gamma$ relative to $\rho$} is a proalgebraic group $\cG$,
defined over $F$, which is an extension of $S$ by a prounipotent
group $\cU$,
$$1\lra\cU\lra\cG\lra S\lra 1,$$
together with a lift $\widetilde{\rho}:\Gamma\ra\cG$ of $\rho$.
The group $\cG$ is required to satisfy the obvious universal
mapping property.  If $S$ is the trivial group then this reduces
to the usual unipotent completion (first defined by Malcev in the
case $F=\zq$).

Examples of this construction may be found in \cite{knudson}.  We
recall some relevant facts that will be used below. Suppose $S$ is
trivial so that we are considering the unipotent $F$-completion.

\begin{prop}\label{filtration} Let $G$ be a group with a
filtration $$G=G^1\supseteq G^2\supseteq G^3\supseteq \cdots$$
satisfying the following conditions.
\begin{enumerate}
\item[$(1)$] The graded quotients $G^i/G^{i+1}$ are
finite-dimensional $F$-vector spaces.

\item[$(2)$] For all $i$, $(G^i/\Gamma^iG)\otimes F=0$.
\end{enumerate}
Then $\cU=\varprojlim G/G^i$ is the unipotent $F$-completion of
$G$.
\end{prop}

\begin{proof}  This is a modification of Proposition 3.5 of
\cite{knudson}.  The unipotent $F$-completion of $G$ is obtained
as the inverse limit of the $F$-completions of each $G/\Gamma^iG$.
Since $(G^i/\Gamma^iG)\otimes F=0$ for all $i$, and since
$G^i/G^{i+1}$ is an $F$-vector space, the completion of
$G/\Gamma^iG$ over $F$ is the group $G/G^i$.  Thus,
$\cU=\varprojlim G/G^i$ is the unipotent $F$-completion of $G$.
\end{proof}

The extension
$$1\lra\cU\lra\cG\lra S\lra 1$$ is split (\cite{knudson}, p.~195).
An obvious question to consider is the relationship between the
group $\cU$, which is prounipotent, and the unipotent completion
of the kernel of $\rho:\Gamma\ra S$.  Denote this kernel by $T$
and let $\cT$ be its unipotent $F$-completion.  Then the map
$T\ra\cU$ induces a map $\Phi:\cT\ra\cU$.  Let $L$ be the image of
$\rho$.

\begin{prop}\label{injective} Suppose that $H_1(T;F)$ is
finite-dimensional.  If the action of $L$ on $H_1(T;F)$ extends to
a rational representation of $S$ $($for example, if $L=S)$, then
the kernel of $\Phi$ is central in $\cT$.
\end{prop}

\begin{proof} See \cite{knudson}, Proposition 4.2.
\end{proof}

\begin{prop}\label{surjective}  Suppose that $H_1(T;F)$ is
finite-dimensional.  If $\rho:\Gamma\ra S$ is surjective, then
$\Phi:\cT\ra\cU$ is surjective.
\end{prop}

\begin{proof} See \cite{knudson}, Proposition 4.3.
\end{proof}

\section{Continuous Cohomology}\label{continuous}
Suppose that $\pi$ is a projective limit of groups,
$$\pi=\varprojlim \pi_\alpha,$$ and let $F$ be a field.  We define
the continuous cohomology of $\pi$ to be
$$H^i_{\text{cts}}(\pi;F) = \varinjlim H^i(\pi_\alpha;F).$$
For example, if $\pi$ is the Galois group of a field extension
$L/K$, then $H^i_{\text{cts}}(\pi;F)$ is simply the usual Galois
cohomology.

This construction is relevant here in light of the following
result.

\begin{prop}\label{h2} Let $\rho:\Gamma\ra S$ be a split
surjective representation and let $\cG$ be the completion relative
to $\rho$.  Assume that $H_1(T;F)$ is finite-dimensional
$(T=\ker\rho)$ and that $\Phi:\cT\ra\cU$ is an isomorphism.  Then
the restriction map $$H^2_{\text{\em cts}}(\cG;F)\lra
H^2(\Gamma;F)$$ is injective.
\end{prop}

\begin{proof} See \cite{knudson}, Corollary 5.5.
\end{proof}

\section{The Completion of $\slna$}\label{slna}
Recall that $k$ denotes a number field and $C$ is a smooth affine
curve over $k$.  The coordinate ring of $C$ is denoted by $A$.
Assume that $x$ is a $k$-rational point of $C$ and let $\mx$ be
the associated maximal ideal of $A$.  We also assume that
$SK_1(A)\otimes \zq=0$.  With these assumptions, we now proceed to
compute the completion of $\slna$, $n\ge 3$.

Let $\rho:\slna\ra\slnk$ be the map induced by reducing modulo
$\mx$.  Note that $\rho$ is split surjective.  Denote by $\ahat$
the $\mx$-adic completion of $A$; this is a complete regular local
ring of dimension $1$ and is thus isomorphic to the power series
ring $k[[T]]$.

\begin{theorem}\label{mainthm} The completion of $\slna$ relative
to $\rho$ is the group $\slnahat$.
\end{theorem}

\begin{proof} Let $K$ be the kernel of $\rho$ and define a central
series $K^\bullet$ by
$$K^i=\{X\in\slna : X\equiv I \mod\mx^i\}.$$  Then for $n\ge 3$,
we have isomorphisms \cite{bass}
$$K^i/E_n(A,\mx^i)\stackrel{\cong}{\lra} SK_1(A,\mx^i)$$
for all $i$.  Here, $E_n$ denotes the subgroup generated by
elementary matrices.  By considering the long exact $K$-theory
sequence associated to $(A,\mx^i)$, our hypothesis that
$SK_1(A)\otimes\zq=0$ yields a surjective map
$$K_2(A/\mx^i)\otimes \zq \lra SK_1(A,\mx^i)\otimes \zq.$$  By
choosing a uniformizing parameter $t$ at $x$, we obtain an
isomorphism of rings $A/\mx^i\cong k[t]/t^i$.  Note that for all
$i$,
\begin{eqnarray*}
K_2(k[t]/t^i) & \cong & K_2(k[t]/t^i,(t))\oplus K_2(k) \\
              & \cong & (\Omega^1_{k|\zz})^{i-1}\oplus K_2(k)
\end{eqnarray*}
(see \cite{graham}) and hence $K_2(k[t]/t^i)\otimes\zq=0$.  Thus,
$SK_1(A,\mx^i)\otimes\zq=0$.

Now, since we have a sequence of inclusions
$$E_n(A,\mx^i)\subseteq \Gamma^iK\subseteq K^i,$$
and since $K^i/E_n(A,\mx^i)$ is torsion, we see that
$K^i/\Gamma^iK$ is torsion.  Note that the graded quotients
$K^i/K^{i+1}$ are finite-dimensional $k$-vector spaces:
$$K^i/K^{i+1}\cong {\mathfrak sl}_n(\mx^i/\mx^{i+1})\cong
{\mathfrak sl}_n(k).$$ By Proposition \ref{filtration}, the
unipotent $k$-completion of $K$ is the group
$$K(\ahat) =\varprojlim K/K^i.$$  This group fits into the split
extension
\begin{equation}\label{seq1}
1\lra K(\ahat)\lra \slnahat\lra \slnk\lra 1.
\end{equation}
Now, if the completion of $\slna$ relative to $\rho$ is the
extension
$$1\lra\cU\lra\cG\lra \slnk\lra 1,$$ we have a homomorphism
$\Phi:K(\ahat)\ra\cU$.  Note that $$H_1(K;k) \cong K/K^2 \oplus
(K/\Gamma^2K \otimes k) \cong K/K^2$$ is a finite-dimensional
$k$-vector space.  By Proposition \ref{surjective} $\Phi$ is
surjective and since the center of $K(\ahat)$ is trivial,
Proposition \ref{injective} implies that $\Phi$ is injective.
Since the extension (\ref{seq1}) is split, we have
$\cG\cong\slnahat$.
\end{proof}

\section{Application to $K$-theory}\label{k2}
For $n\ge 3$, we have the following chain of isomorphisms:
\begin{eqnarray*}
H^2(\slna;k) & \cong & H_2(\slna;k) \\
             & \cong & (H_2(\slna;\zz)\otimes k) \oplus
             \text{Tor}^1(H_1(\slna;\zz),k) \\
             & \cong & (K_2(A)\otimes k) \oplus
             \text{Tor}^1(SK_1(A),k) \\
             & \cong & K_2(A)\otimes k.
\end{eqnarray*}
This is the primary motivation for calculating the completion of
$\slna$.  By Proposition \ref{h2}, we have an injection
$$H^2_{\text{cts}}(\slnahat;k)\lra H^2(\slna;k),$$ and therefore
we obtain a lower bound on $K_2(A)\otimes k$.  Unfortunately, the
lower bound is not useful.

\begin{prop}\label{h2comp} For $n\ge 3$, $H^2_{\text{\em
cts}}(\slnahat;k)=0$.
\end{prop}

\begin{proof} Note that $\ahat\cong k[[T]]$ and so we may as well
consider the curve $C={\mathbb A}^1$, $A=k[t]$.  Then we have the
following:
\begin{eqnarray*}
H^2_{\text{cts}}(SL_n(k[[T]]);k) & \hookrightarrow &
H^2(SL_n(k[t]);k) \\
     & \cong & H_2(SL_n(k[t]);k) \\
     & \cong & H_2(SL_n(k);k) \\
     & \cong & K_2(k)\otimes k \\
     &  = & 0.
\end{eqnarray*}
\end{proof}

\end{document}